\documentclass[letterpaper, 10 pt, conference]{ieeeconf}  

\IEEEoverridecommandlockouts                              
\usepackage{graphicx,amssymb,amstext}
\usepackage{amsmath}
\usepackage{tikz,lscape}
\usepackage{hyperref}
\usepackage{epsfig}
\usetikzlibrary{shapes,arrows}
\usepackage{verbatim}

\title{\LARGE \bf
Robust Estimation of Optical Phase Varying\\
as a Continuous Resonant Process
}

\author{Shibdas Roy$^{1}$*, Ian R. Petersen$^{2}$ and Elanor H. Huntington$^{3}$%
\thanks{$^{1}$S. Roy, $^{2}$I. R. Petersen and $^{3}$E. H. Huntington are with the School of Engineering and Information Technology, University of New South Wales, Canberra.}%
\thanks{*\tt\small shibdas.roy at student.adfa.edu.au}%
}

\begin{document}

\tikzstyle{block} = [draw, fill=white, rectangle, 
    minimum height=1em, minimum width=2em]
\tikzstyle{open}=[inner sep=1mm]
\tikzstyle{none}=[inner sep=0mm]
\tikzstyle{divide} = [draw, fill=white, circle, inner sep=0.5mm]

\maketitle
\thispagestyle{empty}
\pagestyle{empty}

\begin{abstract}

It is well-known that adaptive homodyne estimation of continuously varying optical phase provides superior accuracy in the phase estimate as compared to adaptive or non-adaptive static estimation. However, most phase estimation schemes rely on precise knowledge of the underlying parameters of the system under measurement, and performance deteriorates significantly with changes in these parameters; hence it is desired to develop robust estimation techniques immune to such uncertainties. In related works, we have already shown how adaptive homodyne estimation can be made robust to uncertainty in an underlying parameter of the phase varying as a simplistic Ornstein-Uhlenbeck stochastic noise process. Here, we demonstrate robust phase estimation for a more complicated resonant noise process using a guaranteed cost robust filter.

\end{abstract}

\section{INTRODUCTION}

\bstctlcite{BSTcontrol}

Quantum phase estimation is the problem of estimating an unknown classical phase involved in the dynamics of a quantum system \cite{WM}. Precise phase estimation plays a key role in quantum computation \cite{HWA}, communication \cite{SPK,CHD} and metrology \cite{GLM1}. A fundamental bound on precision is imposed by quantum mechanics \cite{GLM2} and this limits gravitational wave detection \cite{GMM} and can guarantee security in quantum cryptography \cite{IWY}.

Adaptive homodyne \emph{single-shot} measurements of a \emph{fixed} unknown phase can yield mean-square estimation errors below the standard quantum limit (SQL) \cite{HMW,WK1,WK2,MA,BW1}, which is the minimum level of quantum noise obtained using standard measurements without real-time feedback. It is, however, practically more desirable to precisely track a \emph{continuously varying} phase \cite{BW2,TSL,MT,TW,YNW}.

In Refs. \cite{TW,YNW} the signal phase to be estimated is allowed to continuously evolve under the influence of an unmeasured classical stochastic Ornstein-Uhlenbeck (OU) noise process. However, since it is physically unreasonable to precisely know the underlying parameters of the noise process, the estimation process is significantly affected due to the uncertainty in these parameters. Hence, it is desired to make the phase estimation robust to uncertainties in these parameters.

The authors have previously demonstrated in Ref. \cite{RPH1} that using a robust feedback filter designed based on a guaranteed cost robust filtering approach \cite{PM} can improve the estimation process as compared to an optimal filter. In related works, the authors have also shown how continuous phase estimation using smoothing (rather than filtering alone) can be made robust to uncertainties in the phase being measured and/or the measured output for a coherent beam of light \cite{RPH2} and a phase-squeezed beam of light \cite{RPH3}, by employing robust fixed-interval smoothing theory from Ref. \cite{MSP} for continuous uncertain systems admitting a certain integral quadratic constraint.

These works, however, dealt with estimating the signal phase modulated by an Ornstein-Uhlenbeck noise process, which turns out to be a much simplified noise model as compared to the kind of noises that in practice corrupt the signal phase to be estimated. We, therefore, consider a more relevant and complicated second-order resonant noise process here and design a guaranteed cost robust feedback filter based on Ref. \cite{PM} for adaptive homodyne phase estimation of a coherent light beam phase-modulated by the same.

\section{RESONANT NOISE PROCESS}

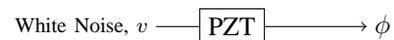
\begin{figure}[!b]
\centering
\begin{tikzpicture}[scale=1]
	\node [open] (white) at (0,0) {\footnotesize White Noise, $v$};
    \node [block] (pzt) at (2,0) {PZT};
    \node [open] (phase) at (4,0) {$\phi$};

    \draw [->] (white) -- (pzt) -- (phase);
\end{tikzpicture}
\caption{Resonant noise}
\label{fig:resonant_noise}
\end{figure}

The resonant noise process we consider here is the one typically generated by a piezo-electric transducer (PZT) driven by an input white noise as in Fig. \ref{fig:resonant_noise}. 

%

\subsection{Transfer Function}\label{sec:resonant_tf}

The simplified transfer function of a typical PZT is given by:
\begin{equation}\label{eq:pzt_tf}
\boxed{G(s) := \frac{\phi}{v} = \frac{\kappa}{s^2+2\zeta\omega_r s+\omega_r^2},}
\end{equation}
where $\kappa$ is the gain, $\zeta$ is the damping factor, $\omega_r$ is the resonant frequency (rad/s) and $v$ is a zero-mean white Gaussian noise with unity amplitude.

\begin{figure}[!t]
\hspace*{-5mm}
\includegraphics[width=0.56\textwidth]{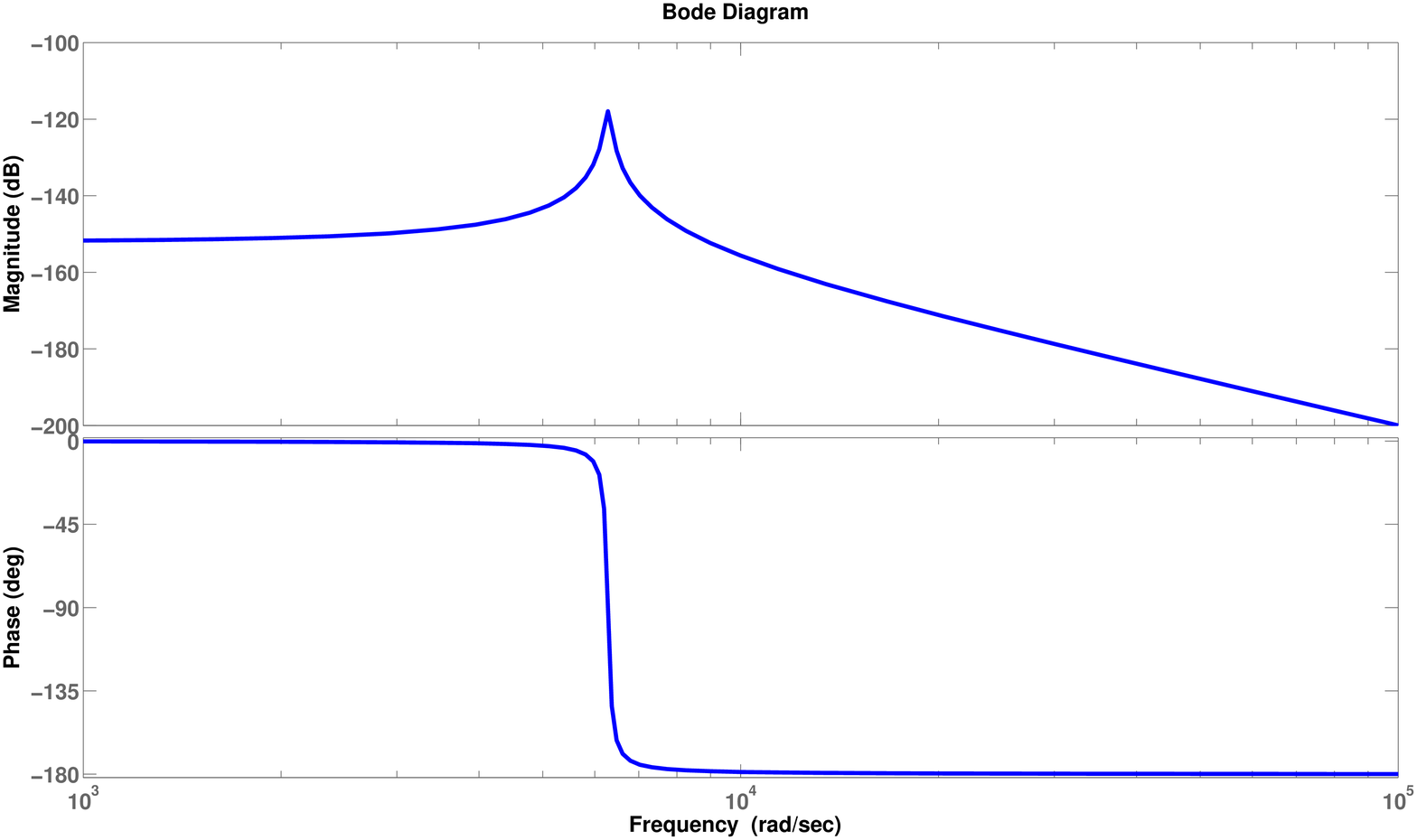}
\caption{Bode plot of the resonant noise process.}
\label{fig:bode_resonant}
\end{figure}

We use the following values for the parameters above: $\kappa = 1$, $\zeta = 0.01$ and $\omega_r = 6.283 \times 10^3$ rad/s, i.e. 1 kHz. The corresponding Bode plot of the transfer function (\ref{eq:pzt_tf}) is given in Fig. \ref{fig:bode_resonant}.

\subsection{State-Space Realization}

From (\ref{eq:pzt_tf}), we obtain:
\begin{equation}
\ddot{\phi}(t)+2\zeta\omega_r\dot{\phi}(t)+\omega_r^2\phi(t) = \kappa v(t).
\end{equation}

Let $x_1 := \phi$ and $x_2 := \dot{\phi}$. Then, we have:
\begin{align}
\dot{x}_1 &= x_2,\\
\dot{x}_2 &= -2\zeta\omega_r x_2 - \omega_r^2 x_1 + \kappa v.
\end{align}

Let $\mathbf{x} := \left[\begin{array}{c}
x_1\\
x_2
\end{array}\right]$.

So, the state-space realization of the PZT output is:
\begin{equation}\label{eq:process_eqn}
\boxed{\mathbf{\dot{x}} = \mathbf{Ax} + \mathbf{G}v,}
\end{equation}
where
\[
\mathbf{A} := \left[\begin{array}{cc}
0 & 1\\
-\omega_r^2 & -2\zeta\omega_r
\end{array}\right], \qquad
\mathbf{G} := \left[\begin{array}{c}
0\\
\kappa
\end{array}\right].
\]

\vspace*{2mm}

\section{OPTIMAL KALMAN FILTER}

The standard adaptive homodyne phase estimation of a coherent state of light is optimal for given values of the parameters when using a Kalman filter in the feedback loop to adapt the phase of the local oscillator. Under a linearization approximation, the homodyne photocurrent from the adaptive phase estimation system is given by \cite{TW}:
\begin{equation}
I(t)dt = 2|\alpha|[\phi(t)-\hat{\phi}(t)]dt + dW(t),
\end{equation}
where $|\alpha|$ is the amplitude of the coherent state with photon flux given by $\mathcal{N}:=|\alpha|^2$, $\hat{\phi}$ is the \emph{intermediate phase estimate}, and $dW$ is Wiener noise arising from quantum vacuum fluctuations.

The \emph{instantaneous estimate} $\theta(t)$ is given by \cite{TW}:
\begin{equation}\label{eq:meas_eqn}
\theta(t) := \hat{\phi}(t)+\frac{I(t)}{2|\alpha|} = \phi(t)+\frac{1}{2|\alpha|}w(t),
\end{equation}
where $w := \frac{dW}{dt}$ is a zero-mean Gaussian white noise with unity amplitude.

Expressing (\ref{eq:meas_eqn}) in terms of $\mathbf{x}$ defined in the previous section, we get the measurement model as:
\begin{equation}
\boxed{\theta = \mathbf{Hx}+\mathbf{J}w,}
\end{equation}
where $\mathbf{H} := \left[\begin{array}{cc} 1 & 0 \end{array}\right]$ and $\mathbf{J} := \left[\frac{1}{2|\alpha|}\right]$.

Rewriting the equations for the system under consideration, we get:
\begin{eqnarray}\label{eq:sys_model}
\boxed{
\begin{split}
\textsf{\small Process model:} \ \ \mathbf{\dot{x}} &= \mathbf{Ax} + \mathbf{G}v, \\
\textsf{\small Measurement model:} \ \ \theta &= \mathbf{Hx} + \mathbf{J}w,
\end{split}}
\end{eqnarray}
where
\begin{align*}
E[v(t)v(\tau)]& = \mathbf{R}\delta(t - \tau), \\
E[w(t)w(\tau)]& = \mathbf{S}\delta(t - \tau), \\
E[v(t)w(\tau)]& = 0.
\end{align*}

Since $v$ and $w$ are unity amplitude white noise processes, both $\mathbf{R}$ and $\mathbf{S}$ are unity (scalars).

The continuous-time algebraic Riccati equation to be solved to construct the steady-state Kalman filter for the system is then \cite{RGB}:
\begin{equation}\label{eq:kalman_riccati}
\mathbf{AP}+\mathbf{PA}^T+\mathbf{GRG}^T-\mathbf{PH}^T(\mathbf{JSJ}^T)^{-1}\mathbf{HP} = \mathbf{0}.
\end{equation}

The Kalman filter equation is \cite{RGB}:
\begin{equation}
\boxed{\mathbf{\dot{\hat{x}}} = (\mathbf{A}-\mathbf{KH})\mathbf{\hat{x}}+\mathbf{KHx}+\mathbf{KJ}w,}
\end{equation}
where $\mathbf{K} := \mathbf{PH}^T(\mathbf{JSJ}^T)^{-1}$ is the Kalman gain.

Using $\kappa = 1$, $\zeta = 0.01$ and $\omega_r = 6.283 \times 10^3$ rad/s, i.e. 1 kHz as in section \ref{sec:resonant_tf} and $|\alpha| = 6 \times 10^8$, we get:
\small
\begin{equation}
\mathbf{P} = \left[\begin{array}{cc}
3.33785970\times10^{-14} & 8.02174133\times10^{-10}\\
8.02174133\times10^{-10} & 3.99751822\times10^{-5}
\end{array}\right],
\end{equation}
\normalsize
and
\begin{equation}
\mathbf{K} = \left[\begin{array}{c}
\small
4.80651797\times10^{4}\\
1.15513075\times10^{9}
\normalsize
\end{array}\right].
\end{equation}

\vspace*{2mm}

\section{ROBUST FILTER}

In this section, we make our filter robust to uncertainties in the parameters underlying the system matrix $\mathbf{A}$ using the guaranteed cost estimation robust filtering approach given in Ref. \cite{PM}.

We introduce uncertainty in $\mathbf{A}$ as follows:
\[
\mathbf{A} \to \mathbf{A} + \left[\begin{array}{cc}
0 & 0\\
-\mu_1\delta_1\omega_r^2 & -2\mu_2\delta_2\zeta\omega_r
\end{array}\right],
\]
where $\Delta := \left[\begin{array}{cc}\delta_1 & \delta_2\end{array}\right]$ is an uncertain parameter satisfying $||\Delta|| \leq 1$, i.e. $\delta_1^2+\delta_2^2 \leq 1$, and $0 \leq \mu_1 < 1$, $0 \leq \mu_2 < 1$ are parameters determining the levels of uncertainty.

The process and measurement models from (\ref{eq:sys_model}) now take the form:
\begin{eqnarray}\label{eq:uncertain_model}
\boxed{
\begin{split}
\textsf{\small Process model:} \ \ \mathbf{\dot{x}} &= (\mathbf{A}+\mathbf{D_1\Delta E_1})\mathbf{x} + \mathbf{G}v, \\
\textsf{\small Measurement model:} \ \ \theta &= \mathbf{Hx} + \mathbf{J}w,
\end{split}}
\end{eqnarray}
where
\[
\mathbf{D_1}:=\left[\begin{array}{c} 0\\ 1\end{array}\right], \qquad \mathbf{E_1}:=\left[\begin{array}{cc}
-\mu_1\omega_r^2 & 0\\
0 & -2\mu_2\zeta\omega_r
\end{array}\right].
\]

As in Ref. \cite{PM}, the Riccati equation to be solved in order to construct the guaranteed cost filter for the system is:
\begin{eqnarray}
\begin{split}
\mathbf{AQ}&+\mathbf{QA}^T+\epsilon\mathbf{QE_1}^T\mathbf{E_1Q}- \epsilon\mathbf{QH}^T(\epsilon\mathbf{JSJ}^T)^{-1}\mathbf{HQ}\\
&+\frac{1}{\epsilon}\mathbf{D_1D_1}^T+\mathbf{GRG}^T=\mathbf{0}.
\end{split}
\end{eqnarray}

\begin{figure}[!b]
\hspace*{-5mm}
\includegraphics[width=0.56\textwidth]{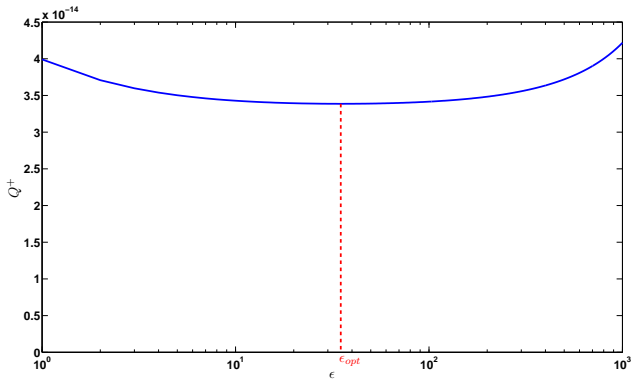}
\caption{$Q^{+}$ versus $\epsilon$ for $\mu_1=0.5$ and $\mu_2=0$.}
\label{fig:optimum_epsilon}
\end{figure}

The stabilising solution of this equation yields an upper bound $\mathbf{\tilde{Q}}$ for the robust filter error covariance. It is desired to obtain the optimum value of $\epsilon$ (that we call $\epsilon_{opt}$) at which the element in the first row and first column of the matrix $\mathbf{\tilde{Q}}$ (that we denote by $Q^{+} := \mathbf{\tilde{Q}}(1,1)$) is minimum. For example, for $\mu_1 = 0.5, \mu_2 = 0$ and nominal values for other parameters, Fig. \ref{fig:optimum_epsilon} shows the plot of $Q^{+}$ versus $\epsilon$, where $\epsilon_{opt}$ is found to be $35$. Thus, we get:

\small
\begin{equation}
\mathbf{\tilde{Q}} = \left[\begin{array}{cc}
3.38608462 \times 10^{-14} & 8.17703018 \times 10^{-10}\\
8.17703018 \times 10^{-10} & 4.09328251 \times 10^{-5}
\end{array}\right].
\end{equation}
\normalsize

The robust filter equation is:
\begin{eqnarray}\label{eq:robust_filter}
\boxed{\begin{split}
\mathbf{\dot{\hat{x}}}&=(\mathbf{A}+\epsilon\mathbf{\tilde{Q}}\mathbf{E_1}^T\mathbf{E_1} - \mathbf{\tilde{Q}}\mathbf{H}^T(\mathbf{JSJ}^T)^{-1}\mathbf{H})\mathbf{\hat{x}}\\
&+\mathbf{\tilde{Q}}\mathbf{H}^T(\mathbf{JSJ}^T)^{-1}\mathbf{Hx} + \mathbf{\tilde{Q}}\mathbf{H}^T(\mathbf{JSJ}^T)^{-1}\mathbf{J}w.
\end{split}}
\end{eqnarray}

\vspace*{2mm}

\section{COMPARISON OF THE ROBUST FILTER WITH THE KALMAN FILTER}

\subsection{Lyapunov Method}\label{sec:lyap_method}

We augment the system given by (\ref{eq:uncertain_model}) with the feedback filter (\ref{eq:robust_filter}) and represent the augmented system by the state-space model:
\begin{equation}\label{eq:ss_model}
\mathbf{\dot{\overline{x}}} = \mathbf{\overline{A}\, \overline{x}} + \mathbf{\overline{B}\, \overline{w}},
\end{equation}
where

\qquad \qquad \( \mathbf{\overline{x}} := 
\left[ \begin{array}{c}
\mathbf{x} \\
\mathbf{\hat{x}}
\end{array} \right] \)
\qquad and \qquad
\( \mathbf{\overline{w}} := 
\left[ \begin{array}{c}
v \\
w
\end{array} \right]. \)

\vspace*{2mm}

Thus, we have:
\[ \mathbf{\overline{A}} = 
\left[ \begin{array}{cc}
\mathbf{A}+\mathbf{D_1\Delta E_1} & \mathbf{0} \\
\mathbf{F} & \mathbf{L}
\end{array} \right], \]
where $\mathbf{F}:=\mathbf{\tilde{Q}}\mathbf{H}^T(\mathbf{JSJ}^T)^{-1}\mathbf{H}$, $\mathbf{L}:=\mathbf{A}+\epsilon\mathbf{\tilde{Q}}\mathbf{E_1}^T\mathbf{E_1} - \mathbf{\tilde{Q}}\mathbf{H}^T(\mathbf{JSJ}^T)^{-1}\mathbf{H}$ and
\[ \mathbf{\overline{B}} = 
\left[ \begin{array}{cc}
\mathbf{G} & \mathbf{0} \\
\mathbf{0} & \mathbf{\tilde{Q}}\mathbf{H}^T(\mathbf{JSJ}^T)^{-1}\mathbf{J}
\end{array} \right]. \]

\vspace*{2mm}

For the continuous-time state-space model (\ref{eq:ss_model}), the steady-state state covariance matrix $\mathbf{P}_S$ is obtained by solving the \emph{Lyapunov equation}:
\begin{equation}\label{eq:lyapunov}
\mathbf{\overline{A}P}_S + \mathbf{P}_S\mathbf{\overline{A}}^T + \mathbf{\overline{B}\, \overline{B}}^T = \mathbf{0},
\end{equation}
where $\mathbf{P}_S$ is the symmetric matrix
\[ 
\mathbf{P}_S := E(\mathbf{\overline{x}\, \overline{x}}^T) :=
\left[ \begin{array}{cc}
\mathbf{P_1} & \mathbf{P_2} \\
\mathbf{P_2}^T & \mathbf{P_3}
\end{array} \right].
\]

The state estimation error can be written as:
\[ \mathbf{e} := \mathbf{x} - \mathbf{\hat{x}} = [\mathbf{1} \, -\mathbf{1}]\mathbf{\overline{x}}, \]
which is mean zero since all of the quantities determining $\mathbf{e}$ are mean zero.

The error covariance matrix is then given as:
\small
\begin{align*} 
\mathbf{\Sigma} :&= E(\mathbf{ee}^T) = [\mathbf{1} \, -\mathbf{1}]E(\mathbf{\overline{x}\, \overline{x}}^T)
\left[ \begin{array}{c}
\mathbf{1} \\
-\mathbf{1}
\end{array} \right] \\ 
&=
[\mathbf{1} \, -\mathbf{1}]
\left[ \begin{array}{cc}
\mathbf{P_1} & \mathbf{P_2} \\
\mathbf{P_2}^T & \mathbf{P_3}
\end{array} \right]
\left[ \begin{array}{c}
\mathbf{1} \\
-\mathbf{1}
\end{array} \right]
= \mathbf{P_1} - \mathbf{P_2} - \mathbf{P_2}^T + \mathbf{P_3}.
\end{align*}
\normalsize

Since we are mainly interested in estimating $x_1 = \phi$, the estimation error covariance of interest is $\sigma^2 = \mathbf{\Sigma}(1,1)$.

\subsection{Standard Quantum Limit}\label{sec:resonant_sql}

The standard quantum limit is set by the minimum error in phase estimation that can be obtained using perfect heterodyne scheme. The SQL for our resonant noise process may be obtained using the same technique employed in Ref. \cite{RPH2} in the case of OU noise process.

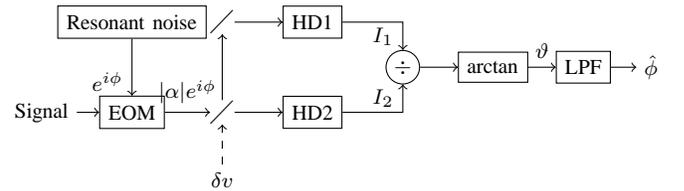
\begin{figure}[!b]
\centering
\begin{tikzpicture}[scale=0.6]
	\node [open] (signal) at (1,0) {\footnotesize Signal};
    \node [block] (eom) at (3,0) {\footnotesize EOM};
	\node [block,align=center] (ou) at (3,2) {\footnotesize Resonant noise};
    \node [none] (bs) at (5,0) {$\diagup$};
    \node [open] (vacuum) at (5,-1.5) {\footnotesize $\delta v$};
    \node [block] (hd2) at (7,0) {\footnotesize HD2};
    \node [none] (mirror) at (5,2) {$\diagup$};
    \node [block] (hd1) at (7,2) {\footnotesize HD1};
    \node [divide] (divide) at (9,1) {\footnotesize $\div$};
    \node [block] (arctan) at (11,1) {\footnotesize arctan};
    \node [block] (lpf) at (13,1) {\footnotesize LPF};
    \node [open] (output) at (14.5,1) {\footnotesize $\hat{\phi}$};
    
    \draw [->] (signal) -- (eom);
    \draw [->] (ou) -- node [left, near end] {\footnotesize $e^{i\phi}$} (eom);
    \draw [->] (eom) -- node [above] {\footnotesize $|\alpha|e^{i\phi}$} (bs);
    \draw [->] (bs) -- (hd2); 
    \draw [->] (hd2) -| node [left, near end] {\footnotesize $I_2$} (divide);
    \draw [->,dashed] (vacuum) -- (bs);
    \draw [->] (bs) -- (mirror);
    \draw [->] (mirror) -- (hd1);
    \draw [->] (hd1) -| node [left, near end] {\footnotesize $I_1$} (divide);
    \draw [->] (divide) -- (arctan);
    \draw [->] (arctan) -- node [above] {\footnotesize $\vartheta$} (lpf);
    \draw [->] (lpf) -- (output);
\end{tikzpicture}
\caption{Block diagram of the dual-homodyne scheme for deducing the SQL for the resonant noise.}
\label{fig:dual_hd_sql}
\end{figure}

We use the fact that the heterodyne scheme of measurement is, in principle, equivalent to, and incurs the same noise penalty as, \emph{dual-homodyne} scheme \cite{TW}, such as the schematic depicted in Fig. \ref{fig:dual_hd_sql}. A signal at the input is phase-modulated using an electro-optic modulator (EOM) that is driven by the resonant noise source. The modulated signal is then split using a $50-50$ beamsplitter into two arms each with a homodyne detector (HD1 and HD2, respectively, with the local oscillator phase of HD1 $\pi/2$ out of phase with that of HD2). The ratio of the output signals of the two arms goes to an \emph{arctan} block, the output of which is fed into a low-pass filter (LPF). The phase estimation error would be minimum when this LPF is an optimal Kalman filter.

The output signals of the two arms are:
\begin{eqnarray*}
I_1 = \frac{1}{\sqrt{2}} \left( 2|\alpha| \sin\phi + n_1 + n_2 \right), \\ 
I_2 = \frac{1}{\sqrt{2}} \left( 2|\alpha| \cos\phi + n_3 - n_4 \right),
\end{eqnarray*}
where $n_1$ and $n_3$ are measurement noises of the two homodyne detectors, respectively, and $n_2$ and $n_4$ are the noises arising from the vacuum entering the empty port of the input beamsplitter corresponding to the two arms, respectively. All these noises are assumed to be zero-mean white Gaussian.

The output of the arctan block is:
\begin{equation}
\vartheta = \arctan \left( \frac{2|\alpha|\sin\phi + n_1 + n_2}{2|\alpha|\cos\phi + n_3 - n_4} \right).
\end{equation}

Assuming the input noises are small, a Taylor series expansion upto first-order terms of the right-hand side yields:
\begin{equation}
\vartheta \approx \phi + \frac{1}{2|\alpha|}n_1 + \frac{1}{2|\alpha|}n_2.
\end{equation}

Expressing this equation in terms of $\mathbf{x}$, we get the measurement model as:
\begin{equation}\label{eq:dual_hd_meas}
\boxed{\vartheta = \mathbf{Hx}+\mathbf{J}w,}
\end{equation}
where $\mathbf{H} = \left[\begin{array}{cc} 1 & 0 \end{array}\right]$ and $\mathbf{J} = \left[\begin{array}{cc} \frac{1}{2|\alpha|} & \frac{1}{2|\alpha|} \end{array}\right]$.

The error covariance matrix of the optimal steady-state Kalman filter for the process given by (\ref{eq:process_eqn}) and the measurement given by (\ref{eq:dual_hd_meas}) may be obtained by solving an algebraic Riccati equation of the form (\ref{eq:kalman_riccati}) for $\mathbf{P}$. The error covariance of interest (i.e. that in estimating $x_1 = \phi$) is then $\sigma^2 = \mathbf{P}(1,1)$.

\subsection{Comparison of Estimation Errors}

\begin{figure}[!b]
\hspace*{-5mm}
\includegraphics[width=0.56\textwidth]{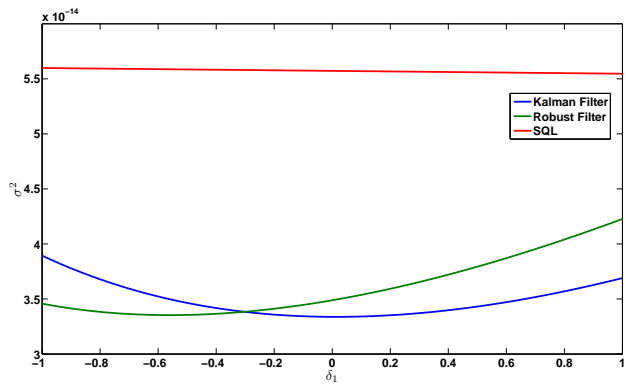}
\caption{Comparison of estimation error covariance $\sigma^2$ for the different filters with $\mu_1=0.2$ and $\mu_2=0$.}
\label{fig:mu1-20}
\end{figure}

\begin{figure}[!htb]
\hspace*{-5mm}
\includegraphics[width=0.56\textwidth]{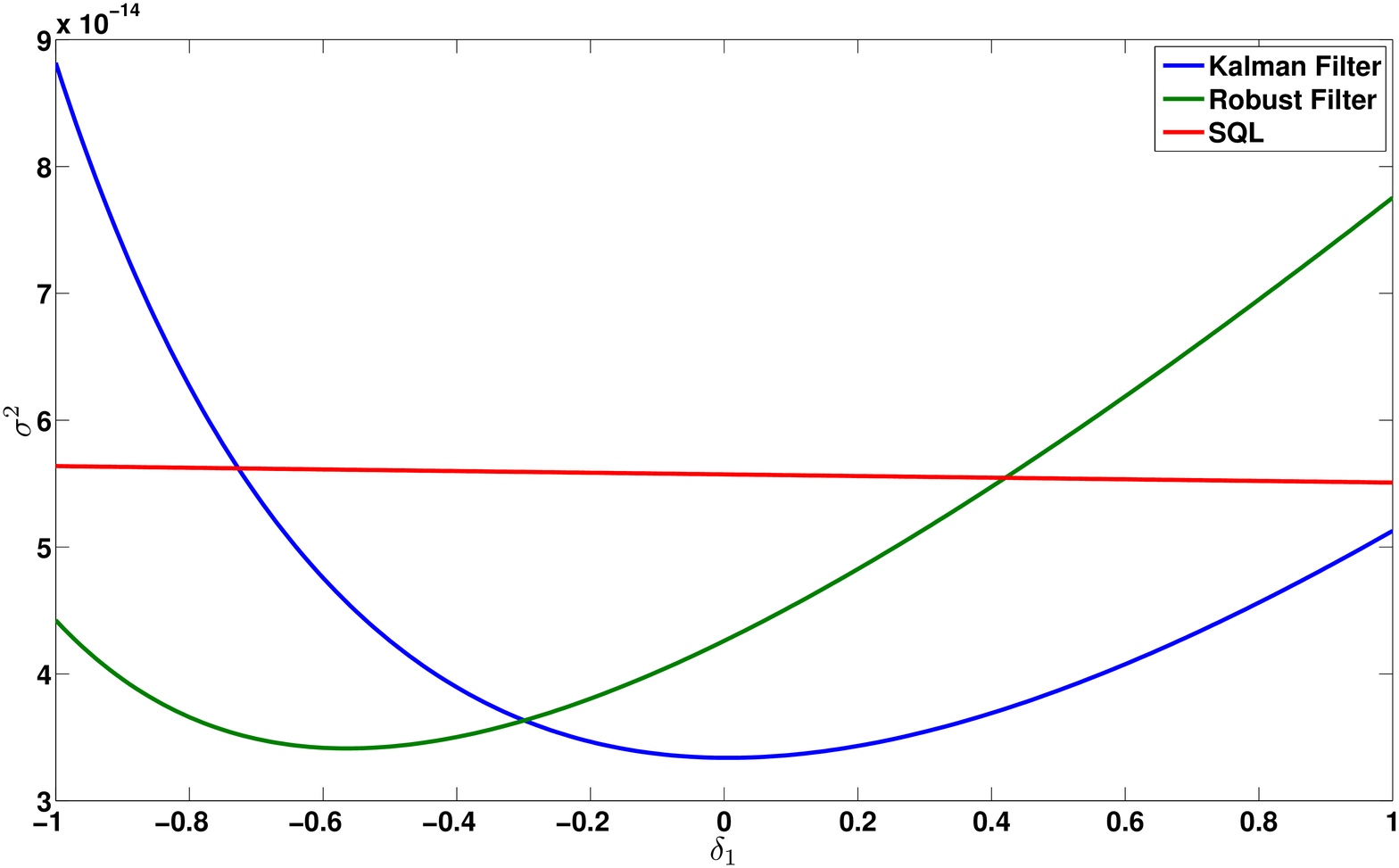}
\caption{Comparison of estimation error covariance $\sigma^2$ for the different filters with $\mu_1=0.5$ and $\mu_2=0$.}
\label{fig:mu1-50}
\end{figure}

\begin{figure}[!htb]
\hspace*{-5mm}
\includegraphics[width=0.56\textwidth]{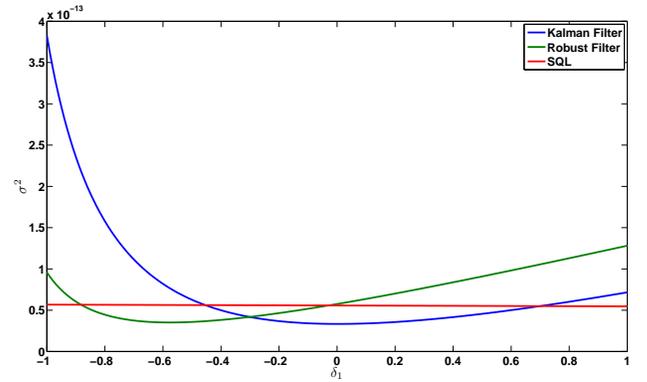}
\caption{Comparison of estimation error covariance $\sigma^2$ for the different filters with $\mu_1=0.8$ and $\mu_2=0$.}
\label{fig:mu1-80}
\end{figure}

The estimation errors may be calculated, as described in section \ref{sec:lyap_method}, for the robust filter, and likewise for the Kalman filter, as a function of $\delta_1$ with $\delta_2=0$ for the nominal values of the parameters and chosen values for $\mu_1$ with $\mu_2=0$. These were used to generate plots of the errors versus $\delta_1$ to be able to compare the performance of the robust filter vis-a-vis the Kalman filter for the uncertain system with respect to the SQL. The SQL value is obtained by designing, as described in section \ref{sec:resonant_sql}, a different Kalman filter for each value of the uncertain parameter. Figs. \ref{fig:mu1-20}, \ref{fig:mu1-50} and \ref{fig:mu1-80} show these plots for $\mu_1=0.2, 0.5$ and $0.8$, respectively.

Clearly, as $\delta_1$ deviates away from $0$ towards $-1$, the performance of the robust filter becomes superior than that of the Kalman filter in relation to the SQL for all levels of $\mu_1$.

\begin{figure}[!b]
\hspace*{-5mm}
\includegraphics[width=0.56\textwidth]{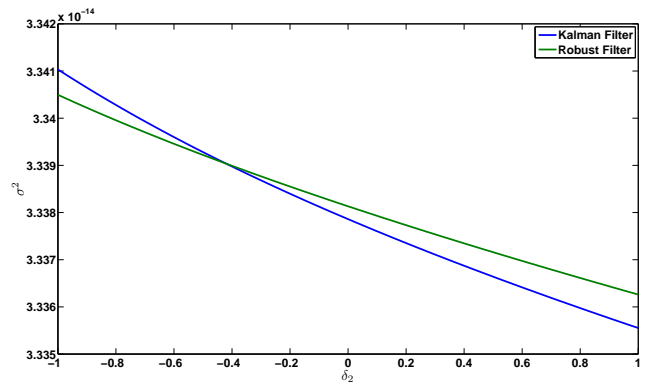}
\caption{Comparison of estimation error covariance $\sigma^2$ for the different filters with $\mu_1=0$ and $\mu_2=0.3$.}
\label{fig:mu2-30}
\end{figure}

\begin{figure}[!htb]
\hspace*{-5mm}
\includegraphics[width=0.56\textwidth]{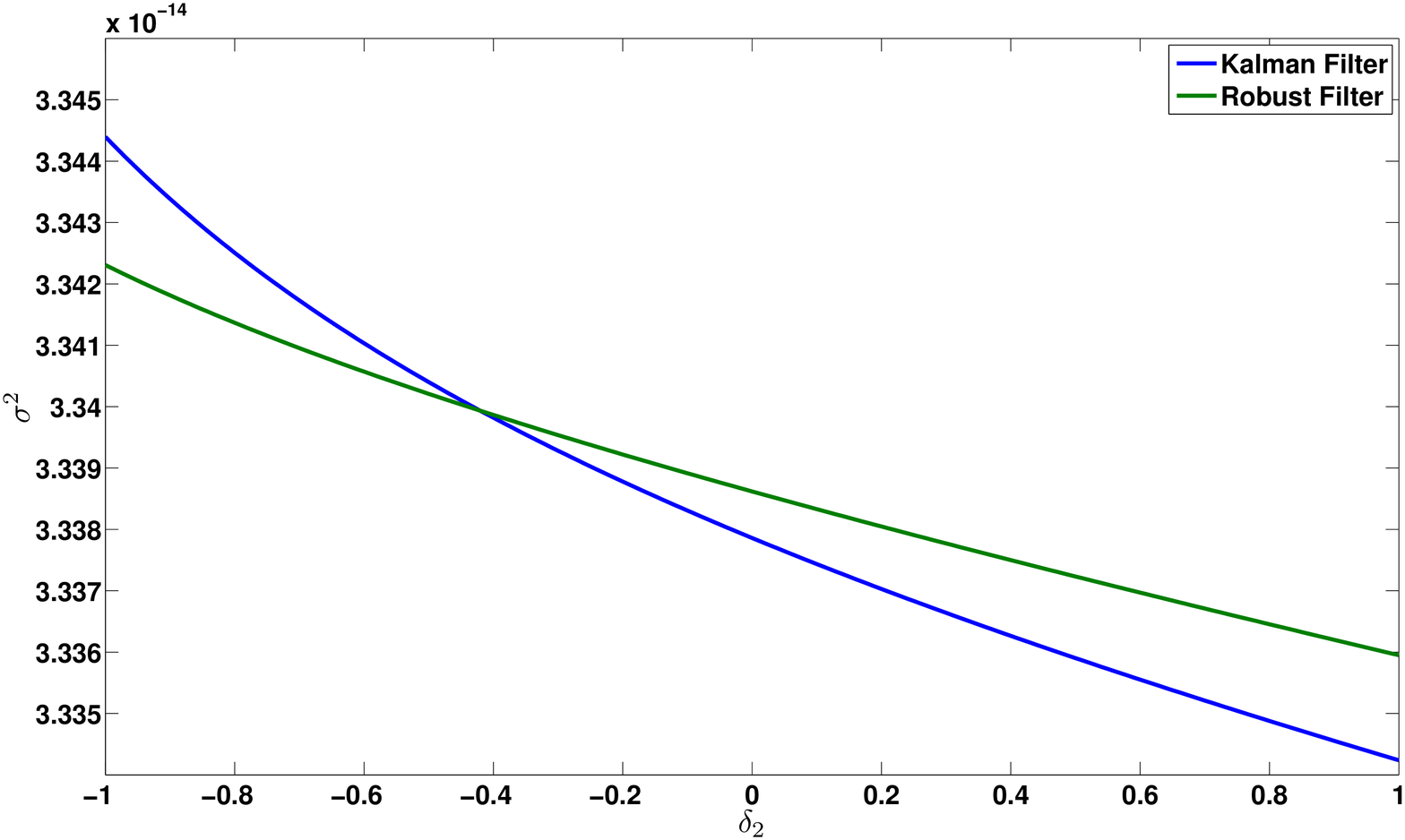}
\caption{Comparison of estimation error covariance $\sigma^2$ for the different filters with $\mu_1=0$ and $\mu_2=0.5$.}
\label{fig:mu2-50}
\end{figure}

\begin{figure}[!htb]
\hspace*{-5mm}
\includegraphics[width=0.56\textwidth]{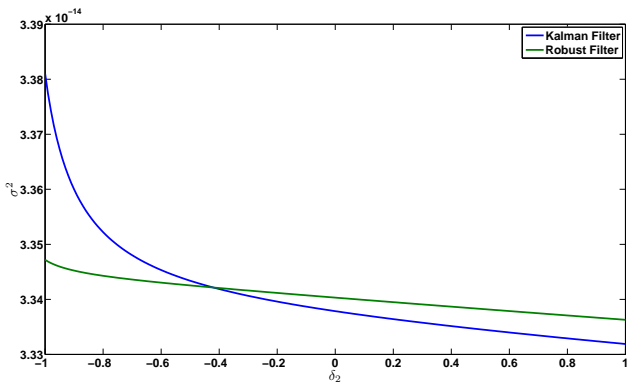}
\caption{Comparison of estimation error covariance $\sigma^2$ for the different filters with $\mu_1=0$ and $\mu_2=0.9$.}
\label{fig:mu2-90}
\end{figure}

Similarly, the estimation errors may be calculated, and plots created, for the robust filter and the Kalman filter as a function of $\delta_2$ with $\delta_1=0$ for the nominal values of the parameters and chosen values for $\mu_2$ with $\mu_1=0$. Figs. \ref{fig:mu2-30}, \ref{fig:mu2-50} and \ref{fig:mu2-90} show these plots for $\mu_2=0.3, 0.5$ and $0.9$, respectively.

In this case too, the robust filter outperforms the Kalman filter as $\delta_2$ tends towards $-1$. The SQL has not been shown in these plots since it is way above these errors.

\vspace*{-2mm}

\section{CONCLUSION}

This paper applies guaranteed cost robust filtering theory to the problem of robustly estimating the optical phase of a coherent state evolving continuously as a classical resonant noise process. We showed the behaviour of the robust filter as compared to the Kalman filter with uncertainties in the important parameters underlying the phase noise. The theory herein may be extended to include robust smoothing rather than filtering alone. Robustness to uncertainties in other parameters, such as the photon flux and noise power, may also be explored. Moreover, robustness for a squeezed state, instead of a coherent state, under the influence of such a resonant noise process would be interesting to investigate.

\bibliographystyle{IEEEtran}
\bibliography{IEEEabrv,bibliography}

\def\url#1{}
\begin{thebibliography}{10}
\providecommand{\url}[1]{#1}
\csname url@samestyle\endcsname
\providecommand{\newblock}{\relax}
\providecommand{\bibinfo}[2]{#2}
\providecommand{\BIBentrySTDinterwordspacing}{\spaceskip=0pt\relax}
\providecommand{\BIBentryALTinterwordstretchfactor}{4}
\providecommand{\BIBentryALTinterwordspacing}{\spaceskip=\fontdimen2\font plus
\BIBentryALTinterwordstretchfactor\fontdimen3\font minus
  \fontdimen4\font\relax}
\providecommand{\BIBforeignlanguage}[2]{{%
\expandafter\ifx\csname l@#1\endcsname\relax
\typeout{** WARNING: IEEEtran.bst: No hyphenation pattern has been}%
\typeout{** loaded for the language `#1'. Using the pattern for}%
\typeout{** the default language instead.}%
\else
\language=\csname l@#1\endcsname
\fi
#2}}
\providecommand{\BIBdecl}{\relax}
\BIBdecl

\bibitem{WM}
H.~M. Wiseman and G.~J. Milburn, \emph{Quantum Measurement and Control}.\hskip
  1em plus 0.5em minus 0.4em\relax Cambridge University Press, 2010.

\bibitem{HWA}
\BIBentryALTinterwordspacing
M.~Hofheinz, H.~Wang, M.~Ansmann, R.~C. Bialczak, E.~Lucero, M.~Neeley, A.~D.
  O'Connell, D.~Sank, J.~Wenner, J.~M. Martinis, and A.~N. Cleland,
  ``Synthesizing arbitrary quantum states in a superconducting resonator,''
  \emph{Nature (London)}, vol. 459, pp. 546--549, March 2009.
  \url{http://www.nature.com/nature/journal/v459/n7246/full/nature08005.html}
\BIBentrySTDinterwordspacing

\bibitem{SPK}
\BIBentryALTinterwordspacing
R.~Slavik, F.~Parmigiani, J.~Kakande, C.~Lundstrom, M.~Sjodin, P.~A. Andrekson,
  R.~Weerasuriya, S.~Sygletos, A.~D. Ellis, L.~Gruner-Nielsen, D.~Jakobsen,
  S.~Herstrom, R.~Phelan, J.~O'Gorman, A.~Bogris, D.~Syvridis, S.~Dasgupta,
  P.~Petropoulos, and D.~J. Richardson, ``All-optical phase and amplitude
  regenerator for next-generation telecommunications systems,'' \emph{Nature
  Photonics}, vol.~4, pp. 690--695, September 2010.
  \url{http://www.nature.com/nphoton/journal/v4/n10/full/nphoton.2010.203.html}
\BIBentrySTDinterwordspacing

\bibitem{CHD}
\BIBentryALTinterwordspacing
J.~Chen, J.~L. Habif, Z.~Dutton, R.~Lazarus, and S.~Guha, ``Optical codeword
  demodulation with error rates below the standard quantum limit using a
  conditional nulling receiver,'' \emph{Nature Photonics}, vol.~6, p. 374, May
  2012.
  \url{http://www.nature.com/nphoton/journal/v6/n6/full/nphoton.2012.113.html}
\BIBentrySTDinterwordspacing

\bibitem{GLM1}
\BIBentryALTinterwordspacing
V.~Giovannetti, S.~Lloyd, and L.~Maccone, ``Advances in quantum metrology,''
  \emph{Nature Photonics}, vol.~5, p. 222, March 2011.
  \url{http://www.nature.com/nphoton/journal/v5/n4/full/nphoton.2011.35.html}
\BIBentrySTDinterwordspacing

\bibitem{GLM2}
\BIBentryALTinterwordspacing
V.~Giovannetti, S.~Lloyd, and L.~Maccone, ``Quantum-enhanced measurements:
  Beating the standard quantum limit,'' \emph{Science}, vol. 306, no. 5700, pp.
  1330--1336, November 2004.
  \url{http://www.sciencemag.org/content/306/5700/1330.abstract}
\BIBentrySTDinterwordspacing

\bibitem{GMM}
\BIBentryALTinterwordspacing
K.~Goda, O.~Miyakawa, E.~E. Mikhailov, S.~Saraf, R.~Adhikari, K.~McKenzie,
  R.~Ward, S.~Vass, A.~J. Weinstein, and N.~Mavalvala, ``A quantum-enhanced
  prototype gravitational-wave detector,'' \emph{Nature Physics}, vol.~4, pp.
  472--476, March 2008.
  \url{http://www.nature.com/nphys/journal/v4/n6/abs/nphys920.html}
\BIBentrySTDinterwordspacing

\bibitem{IWY}
\BIBentryALTinterwordspacing
K.~Inoue, E.~Waks, and Y.~Yamamoto, ``Differential phase shift quantum key
  distribution,'' \emph{Physical Review Letters}, vol.~89, p. 037902, June
  2002.  \url{http://link.aps.org/doi/10.1103/PhysRevLett.89.037902}
\BIBentrySTDinterwordspacing

\bibitem{HMW}
\BIBentryALTinterwordspacing
H.~M. Wiseman, ``Adaptive phase measurements of optical modes: Going beyond the
  marginal q distribution,'' \emph{Physical Review Letters}, vol.~75, pp.
  4587--4590, December 1995.
  \url{http://link.aps.org/doi/10.1103/PhysRevLett.75.4587}
\BIBentrySTDinterwordspacing

\bibitem{WK1}
\BIBentryALTinterwordspacing
H.~M. Wiseman and R.~B. Killip, ``Adaptive single-shot phase measurements: A
  semiclassical approach,'' \emph{Physical Review A}, vol.~56, pp. 944--957,
  July 1997.  \url{http://link.aps.org/doi/10.1103/PhysRevA.56.944}
\BIBentrySTDinterwordspacing

\bibitem{WK2}
\BIBentryALTinterwordspacing
H.~M. Wiseman and R.~B. Killip, ``Adaptive single-shot phase measurements: The
  full quantum theory,'' \emph{Physical Review A}, vol.~57, pp. 2169--2185,
  March 1998.  \url{http://link.aps.org/doi/10.1103/PhysRevA.57.2169}
\BIBentrySTDinterwordspacing

\bibitem{MA}
\BIBentryALTinterwordspacing
M.~A. Armen, J.~K. Au, J.~K. Stockton, A.~C. Doherty, and H.~Mabuchi,
  ``Adaptive homodyne measurement of optical phase,'' \emph{Physical Review
  Letters}, vol.~89, p. 133602, September 2002.
  \url{http://link.aps.org/doi/10.1103/PhysRevLett.89.133602}
\BIBentrySTDinterwordspacing

\bibitem{BW1}
\BIBentryALTinterwordspacing
D.~W. Berry and H.~M. Wiseman, ``Phase measurements at the theoretical limit,''
  \emph{Physical Review A}, vol.~63, p. 013813, December 2000.
  \url{http://link.aps.org/doi/10.1103/PhysRevA.63.013813}
\BIBentrySTDinterwordspacing

\bibitem{BW2}
\BIBentryALTinterwordspacing
D.~W. Berry and H.~M. Wiseman, ``Adaptive quantum measurements of a
  continuously varying phase,'' \emph{Physical Review A}, vol.~65, p. 043803,
  March 2002.  \url{http://link.aps.org/doi/10.1103/PhysRevA.65.043803}
\BIBentrySTDinterwordspacing

\bibitem{TSL}
\BIBentryALTinterwordspacing
M.~Tsang, J.~H. Shapiro, and S.~Lloyd, ``Quantum theory of optical temporal
  phase and instantaneous frequency. ii. continuous-time limit and
  state-variable approach to phase-locked loop design,'' \emph{Physical Review
  A}, vol.~79, p. 053843, May 2009.
  \url{http://link.aps.org/doi/10.1103/PhysRevA.79.053843}
\BIBentrySTDinterwordspacing

\bibitem{MT}
\BIBentryALTinterwordspacing
M.~Tsang, ``Time-symmetric quantum theory of smoothing,'' \emph{Physical Review
  Letters}, vol. 102, p. 250403, June 2009.
  \url{http://link.aps.org/doi/10.1103/PhysRevLett.102.250403}
\BIBentrySTDinterwordspacing

\bibitem{TW}
\BIBentryALTinterwordspacing
T.~A. Wheatley, D.~W. Berry, H.~Yonezawa, D.~Nakane, H.~Arao, D.~T. Pope, T.~C.
  Ralph, H.~M. Wiseman, A.~Furusawa, and E.~H. Huntington, ``Adaptive optical
  phase estimation using time-symmetric quantum smoothing,'' \emph{Physical
  Review Letters}, vol. 104, p. 093601, March 2010.
  \url{http://link.aps.org/doi/10.1103/PhysRevLett.104.093601}
\BIBentrySTDinterwordspacing

\bibitem{YNW}
\BIBentryALTinterwordspacing
H.~Yonezawa, D.~Nakane, T.~A. Wheatley, K.~Iwasawa, S.~Takeda, H.~Arao,
  K.~Ohki, K.~Tsumura, D.~W. Berry, T.~C. Ralph, H.~M. Wiseman, E.~H.
  Huntington, and A.~Furusawa, ``Quantum-enhanced optical-phase tracking,''
  \emph{Science}, vol. 337, no. 6101, p. 1514, September 2012.
  \url{http://www.sciencemag.org/content/337/6101/1514.abstract}
\BIBentrySTDinterwordspacing

\bibitem{RPH1}
S.~Roy, I.~R. Petersen, and E.~H. Huntington, ``Robust filtering for adaptive
  homodyne estimation of continuously varying optical phase,''
  \emph{Proceedings of the 2012 Australian Control Conference}, pp. 454--458,
  November 2012.

\bibitem{PM}
\BIBentryALTinterwordspacing
I.~R. Petersen and D.~C. McFarlane, ``Optimal guaranteed cost control and
  filtering for uncertain linear systems,'' \emph{IEEE Transactions on
  Automatic Control}, vol.~39, no.~9, pp. 1971--1977, September 1994.
  \url{http://ieeexplore.ieee.org/iel4/9/7643/00317138.pdf?arnumber=317138}
\BIBentrySTDinterwordspacing

\bibitem{RPH2}
S.~Roy, I.~R. Petersen, and E.~H. Huntington, ``Adaptive continuous homodyne
  phase estimation using robust fixed-interval smoothing,'' \emph{To appear in
  the Proceedings of the American Control Conference, 2013}.

\bibitem{RPH3}
S.~Roy, I.~R. Petersen, and E.~H. Huntington, ``Robust phase estimation of
  squeezed state,'' \emph{To appear in the Proceedings of the CLEO:QELS 2013}.

\bibitem{MSP}
\BIBentryALTinterwordspacing
S.~O.~R. Moheimani, A.~V. Savkin, and I.~R. Petersen, ``Robust filtering,
  prediction, smoothing, and observability of uncertain systems,'' \emph{IEEE
  Trans. on Circuits and Systems I - Fundamental Theory and Appl.}, vol.~45,
  no.~4, p. 446, April 1998.  \url{http://ieeexplore.ieee.org}
\BIBentrySTDinterwordspacing

\bibitem{RGB}
R.~G. Brown, \emph{Introduction to Random Signal Analysis and Kalman
  Filtering}.\hskip 1em plus 0.5em minus 0.4em\relax John Wiley \& Sons, 1983.

\end{thebibliography}

\end{document}